\documentclass{article}

\usepackage{arxiv}

\usepackage[utf8]{inputenc} 
\usepackage[T1]{fontenc}    
\usepackage{hyperref}       
\usepackage{url}            
\usepackage{booktabs}       
\usepackage{amsfonts}       
\usepackage{nicefrac}       
\usepackage{microtype}      
\usepackage{lipsum}		
\usepackage{graphicx}
\usepackage{natbib}
\usepackage{doi}
\usepackage{bm}
\usepackage{amsthm}
\usepackage{amsmath}
\usepackage{subcaption}
\usepackage{mathtools}
\usepackage{multirow}

\newtheorem{assumption}{Assumption}[section]

\title{OPO: Making Decision-Focused Data Acquisition Decisions}

\date{} 					

\author{ \href{https://orcid.org/0000-0003-3432-2123}{
\includegraphics[scale=0.06]{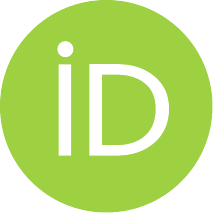}\hspace{1mm}Egon Per\v{s}ak}\\
School of Mathematics \\
University of Edinburgh \\
\texttt{E.Persak@sms.ed.ac.uk} \\
\And
\href{https://orcid.org/0000-0002-8258-9116}{
\includegraphics[scale=0.06]{orcid.pdf}\hspace{1mm}Miguel F. Anjos}\\
School of Mathematics \\
University of Edinburgh \\}

\renewcommand{\shorttitle}{Draft}

\begin{document}
\maketitle

\begin{abstract}
We propose a model for making data acquisition decisions for variables in contextual stochastic optimisation problems. Data acquisition decisions are typically treated as separate and fixed. We explore problem settings in which the acquisition of contextual variables is costly and consequently constrained. The data acquisition problem is often solved heuristically for proxy objectives such as coverage. The more intuitive objective is the downstream decision quality as a result of data acquisition decisions. The whole pipeline can be characterised as an optimise-then-predict-then-optimise (OPO) problem. Analogously, much recent research has focused on how to integrate prediction and optimisation (PO) in the form of decision-focused learning. We propose leveraging differentiable optimisation to extend the integration to data acquisition. We solve the data acquisition problem with well-defined constraints by learning a surrogate linear objective function. We demonstrate an application of this model on a shortest path problem for which we first have to set a drone reconnaissance strategy to capture image segments serving as inputs to a model that predicts travel costs. We ablate the problem with a number of training modalities and demonstrate that the differentiable optimisation approach outperforms random search strategies.
\end{abstract}

\section{Introduction}
Deciding what data to acquire is an integral step in decision making pipelines. Data acquisition (DA) often comes with constraints on how the data can be collected and how costly that collection is in terms of limited resources. We focus on data acquisition in terms of contextual variables which are then used to predict the uncertain parameters for a given contextual decision problem. The actual parameters are observed after the decision is taken. The DA problem as a stand-alone problem does not have an explicit objective in contrast to downstream prediction and decision optimisation problems. A consequence of this ambiguity is that the DA problem is often solved heuristically based on intuition and empirical experience. Not considering the end-to-end effects of DA carries a risk of objective misalignment which potentially limits the capabilities of downstream models.

A more natural setting is to optimise the DA for downstream prediction or decision objectives. This is akin to determining the relative value or shadow price of having access to contextual information. This paper focuses on the relative value case, but the methodology is extendable to the shadow price case. Relative value is relevant in situations where we are trying to optimise a downstream outcome given a fixed and limiting set of DA resources. The idea of shadow prices is relevant when we are trying to decide what costly DA resources to acquire. There are obviously eventual diminishing returns to adding more DA capacity. It is difficult if not impossible to determine the form of returns on additional sensors analytically in a data-driven context. This ambiguity of value often results in suboptimal data acquisition strategies and consequently a reluctance to adopting complex monitoring systems. 

To solve the decision alignment problem for data acquisition we propose OPO, a model which leverages differentiable optimisation (DO) to learn a surrogate data acquisition problem. OPO stands for optimise-predict-optimise, a three-stage modelling approach which segments a data-driven decision-making approach into three stages: data acquisition as an optimisation problem, predicting unknown decision task parameters, and solving the downstream decision task. In this work we explore permutations of melding these tasks including end-to-end training, optimisation and random search strategies over pretrained prediction models, and warm-start strategies for the DA problem.

Our main contributions are:
\begin{itemize}
    \item A fully differentiable model with two non-sequential implicit layers which extends the end-to-end decision-making training to the DA stage. OPO is capable of approximating good solutions to constrained data acquisition problems and is the first decision-focused model for data acquisition.
    \item An implementation of the model on for the task of having to decide a drone reconnaissance strategy to provide inputs for a computer vision model which predicts the travel cost parameters to a downstream shortest path problem.
    \item Empirically evaluating a number of ablations and heuristic search strategies to compare the relative performance of varying degrees of end-to-end training. The best performance is achieved by a combination of warm-starting and finetuning an end-to-end optimisation approach.
\end{itemize}
\section{Optimise-Predict-Optimise Problems}
The problem naturally separates into three stages: data acquisition (DA), prediction (P), and contextual decision task (CDT). Each stage has its own set of decision variables: DA-data acquisition decisions $\textbf{s}$, P-prediction model parameters $W$ including imputation parameters $\Lambda \subset W$, and CDT-decisions in contextual problem $\textbf{x}$. The decision variables in stages DA and P are constant, whereas The decisions in CDT are adjusted for the observed context $\textbf{z}(\textbf{s})$ based on a prediction $\hat{\theta}$ of the true parameters $\theta$.

\subsection{Modelling DA}
The mapping of DA decisions to downstream objectives is implicit and depends on the outcome of optimisation processes. If the objective is prediction quality the mapping is implicit in the sense that the prediction model has to be optimised. If the objective is decision quality then the mapping is the result of two optimisation processes. An exact formulation of this mapping is unknown since we do not have access to the true conditional distribution $\theta |z$. To optimise DA decision-making we need to define some process by which a set of learnable parameters is mapped to DA decisions in accordance with constraints on DA. We assume the constraints are known a priori. Denote the data acquisition decision variables as $\textbf{s}$ where $\textbf{s} \in \{0,1\}^k$. Each $s_i$ is associated with a set of contextual variables contextual variables $\textbf{z}_i \in \mathcal{R}^{c_i}$ and the total number of contextual variables is $\sum_{i=1}^k c_i = C$. Denote the constraints on data acquisition in generic form as $A(\textbf{s}) \leq 0$. An example of such constraints are routing constraints if the DA has a geographic component. Let $p(\textbf{s})$ be the unknown function which maps the DA decision to the expected value of the downstream value either in terms of decision or prediction quality. The true DA problem is then:
\begin{equation}
\label{abstract DA}
\begin{array}{cc}
    \mathrm{arg}\max & p(\textbf{s}) \\
     s.t. & A(\textbf{s}) \leq 0 \: .
\end{array}
\end{equation}
This problem is not computable as the mapping $p(\cdot)$ is unknown. We propose the use of a surrogate linear objective parameterised by vector $\bm{\pi}$ to approximate the optimal solution to the true problem:
\begin{equation}
\label{surrogate DA}
\begin{array}{cc}
    \mathrm{arg}\max & \langle \bm{\pi}, \textbf{s} \rangle \\
     s.t. & A(\textbf{s}) \leq 0 \: .
\end{array}
\end{equation}
If there is any interaction between the variables in $\textbf{s}$ then $p$ is not linear. However, the linear surrogate can obtain the same optimal solution as the original problem so long as the optimal solution to the original problem is on the boundary of the convex hull of the feasible set of the original problem. We believe that this is a reasonable assumption to make.

We operate in a setting where we have access to past observations $(z_i,\theta_i)$ which we use to optimise the surrogate cost $\bm{\pi}$. Given an imputation and prediction strategy, OPO computes the sample downstream objective value as a proxy for $p(\cdot)$. The solution to the surrogate problem \ref{surrogate DA} is fed into a prediction model with an imputation process, and the prediction is then used to solve the CDT. We can optimise this system end-to-end using differentiable optimisation. The surrogate cost $\bm{\pi}$ can be interpreted as learning the relative unobserved value of acquiring contextual information. 

\paragraph{Contrast with Feature Selection}
Selecting what subset of variables to observe for a given model can be interpreted as a form of feature selection. Feature selection is often used in classical machine learning to reduce overfitting. It is solved with a range of heuristics scored by empirical evaluations, and at most contains trivial constraints such as an upper bound on the number of included variables. In contrast OPO can operate with any solver-amenable constraints and uses gradient descent instead of heuristics.

\subsection{Prediction}
The notation in this section is specific to using a vision transformer (ViT) as the prediction model for a computer vision prediction task but OPO works for any differentiable prediction model with a learnable imputation strategy.

Let $\textbf{s}(\bm{\pi}) = \textbf{s}^* \in \{0,1\}^k$ be the solution to problem (\ref{surrogate DA}), $\textbf{z} \in \mathcal{R}^M$ a vector of all corresponding contextual variables. The next step is to conditionally predict $p$ unknown parameter values $\theta \in \mathcal{R}^p$. Let $m(\textbf{z}, \textbf{s}^*;W) = \hat{\theta}$ be a prediction model such as a neural network parameterised by $W$. 

We construct a prediction model which is designed to work with masked inputs. This is mainly to induce stability in prediction and enable the same model to be reasonably usable for any DA solution. The core of the prediction model is the encoder used by \citet{xie2022simmim}. Each segment of the image is tokenised using a convolutional layer with $d_{m}$ kernels where $d_{m}$ is the embedding size. The stride and kernel size are set such that each value in each segment is passed only once. Under this condition the convolutional layer is equivalent to a linear layer. In addition the encoder keeps two more learnable parameters: a class token, and a mask token which we denote as $\Lambda$: both are vectors of size $d_m$. The class token and is always used as the first token providing an anchor for which learning always happens regardless of the selection of tokens. The mask token replaces any tokens corresponding to segments which are not observed. Denote the tokenization layer as $t(\textbf{z}): \mathcal{R}^{c \times k} \mapsto \mathcal{R}^{k \times d_m}$. In the original implementation this operation is not differentiable as it uses an indexing approach. Our equivalent differentiable masked tokenization operation is:
\begin{equation}
\label{eq:masking}
I(\textbf{z},\textbf{s}^*;\Lambda) = \operatorname{diag}(\textbf{s}^*)t(\textbf{z}) + (\bm{1}-\textbf{s}^*)\otimes \Lambda \:.
\end{equation}
The left term copies all tokens which correspond to observed variables, the right terms copies the mask token into all remaining positions. After prepending the class token, the tokens are then passed into the standard encoder module including positional encoding, multi-head attention, skip-connections, fully-connected layers, and layer normalisation. After three transformer blocks we first apply a linear layer to mix across the token dimension and compress the size so the final layer does not present as much of a bottleneck. After a ReLU activation function, we flatten the compressed tokens and latent dimension and pass it to a single linear layer which projects to the size of the unknown parameter vector $\mathcal{R}^p$. Denote the model computation after the masked tokenisation operation as $m_E(\cdot)$. The whole prediction model can therefore be expressed as: 
\begin{equation}
    m(\textbf{z}, \textbf{s}^*;W) =m_E(I(\textbf{z},\textbf{s}^*;\Lambda);W) \:.
\end{equation}

OPO can be seen as a method for constrained (vision) token selection. We leave it to future work to examine how this could be extended to language tasks, but given that reasoning can be expressed as a combinatorial optimisation problem we hypothesise that it is possible.

\subsection{Contextual Decision Task}
Denote the decision variables for that task as $\textbf{x} \in \mathcal{R}^p$. We assume the constraints $B(\textbf{x})$ are certain and that the uncertainty is only in the objective and that the objective is linear. Note that this problem is contextual and changes from instance to instance. Stage CDT is defined as:
\begin{equation}
\label{stage III problem}
\begin{array}{ccc}
    \textbf{x}^*(\theta) = &\mathrm{arg}\min & \langle \theta, \textbf{x} \rangle \\
     &s.t. & B(\textbf{x}) \leq 0 \:.
\end{array}
\end{equation}

\subsection{DFL or PFL}
The Decision-Focused Learning (DFL) approach trains a prediction model to optimise for some measure of downstream decision loss \citep{mandi2023decision}. It has gained a lot of attention as a method for solving contextual optimisation problems. DFL performance in messy settings such as missing data or masked inputs is poorly understood \citep{johnson2023modeling}. Differentiable robust optimisation-type methods have been proposed to deal with label drift/noise \citep{johnson2023characterizing,pervsak2024learning}. To the best of our knowledge, no work has studied the effect of missing values on DFL performance. In contrast, prediction-focused learning (PFL) trains a prediction model to minimise the expected value of some measure of distance between predicted parameters $\hat{\theta}$ and true parameters $\theta$. It has been more broadly studied and entire classes of models have emerged to deal with masked inputs such as the encoder in our prediction model. Given the lack of understanding of the effect of masked inputs on DFL, a natural ablation is to test both training modalities and their performance on the CDT.

\paragraph{PFL Loss} We use mean squared error (MSE) as our PFL loss measure. The loss on a sample of size $n$ is:
\begin{equation}
\label{PFL loss}
    \mathrm{loss_{PFL}}(\theta,\hat{\theta}) = \frac{1}{n}\sum_{i=1}^n (\theta_i - \hat{\theta}_i)^2
\end{equation}

\paragraph{DFL Loss} We use empirical regret as a measure of performance which is the loss in objective value from decisions obtained based on predicted parameters given true parameters.  Many techniques have been developed to calculate or approximate gradients across a constrained optimisation operator $\textbf{x}*(\cdot)$; we discuss specifics in Section \ref{diff}. 

\begin{equation}
\label{DFL loss}
    \mathrm{loss_{DFL}}(\theta,\hat{\theta}) = \frac{1}{n}\sum_{i=1}^n \langle \theta_i, \textbf{x}^*(\hat{\theta}_i)  \rangle - \langle \theta_i, \textbf{x}^*(\theta_i)  \rangle
\end{equation}

\section{Optimising the Surrogate}
\subsection{Differentiation}
\label{diff}
The whole OPO pipeline from DA surrogate parameters $\bm{\pi}$, prediction parameters $\Lambda,W$, and contextual information $\textbf{z}$ to contextual decisions $\textbf{x}$ expressed functionally is
\begin{equation}
    \textbf{x}^*(\textbf{z};\bm{\pi},\Lambda,W) = \textbf{x}^*(
    m_E(I(\textbf{z},
    \textbf{s}(\bm{\pi});
    \Lambda);
    W)
    )
\end{equation}
We can compute the gradient of DFL loss with respect to $\bm{\pi}$ using the chain rule to break down the gradient computation:
\begin{equation}
\label{all gradients}
\frac{\partial \mathrm{loss_{DFL}}}{\partial \bm{\pi}}
     = \frac{\partial \langle \theta, \textbf{x}^*  \rangle}{\partial \textbf{x}^*}
     \underbrace{\frac{\partial \textbf{x}^*}{\partial \hat{\theta}}}_{\text{DO}}
     {\frac{\partial m_E}{\partial I}}
     {\frac{\partial I}{\partial \textbf{s}^*}}
     \underbrace{\frac{\partial \textbf{s}^*}{\partial \bm{\pi}}}_{\text{DO}} \:.
\end{equation}
The first and fourth terms can be obtained using standard analytical differentiation/backpropagation.
\paragraph{Counterfactual Gradient}
The gradient $\frac{\partial I}{\partial \textbf{s}^*}$ needs closer examination as its computation requires knowledge of the whole $\textbf{z}$, not just the subset we decide to observe. We make the following assumption:
\begin{assumption}[Full Observation Assumption]
    The training set consists of fully observed contextual vectors $\textbf{z}$ and fully observed true parameters $\theta$. 
\end{assumption}
We believe this is a reasonable assumption in many settings. For example if a government has engaged in a census-level effort and a decision maker wants to build prediction models based on census-level data but does not have the resources to monitor as many inputs. The assumption is necessary during a backwards pass through OPO. In the forward pass the operations in $I(\cdot)$ which involve $\textbf{z}$ can simply be replaced by an already masked vector which reflects what has actually been observed. In the backwards pass the relevant gradient is $\frac{\partial I}{\partial \tilde{\textbf{s}}}$ which we name the counterfactual gradient. Denote the outputs of $I(\cdot)$ as $\left[t_1 \cdot\cdot\cdot t_k\right]^T$ where $t_i \in \mathcal{R}^{d_m}$ is the token corresponding to DA decision $\textbf{s}_i$. We observe that $I$ is linear in $\textbf{s}$ and that any individual $\textbf{s}_i$ only affects its corresponding token $t_i$. The counterfactual gradient is calculated as:
\begin{equation}
\label{counterfactual gradient}
\begin{array}{c}
\frac{\partial I}{\partial \tilde{\textbf{s}}} = \left[ \nabla t_1 \cdot\cdot\cdot \nabla t_k \right]^T \\ \frac{\partial t_i}{\partial \textbf{s}_j} =
\begin{cases}
    \frac{\partial}{\partial \textbf{s}_j}  \left[ 
\textbf{s}_i t(\textbf{z})_i + (1-\textbf{s}_i)\Lambda \right] = t(\textbf{z})_i - \Lambda & \text{if $i=j$}\\
0 & \text{otherwise} \\
\end{cases}
\end{array}
\end{equation}
We see that this gradient does not just depend on what was actually observed. In the backward pass the imputation layer needs the counterfactual to compute how the prediction would change as a result of a change in data acquisition.

\paragraph{Differentiable Optimisation}
Both terms marked as DO in equation (\ref{all gradients}) require finding a gradient from an optimisation problem with a linear objective. The constraints $A(\textbf{s)}$ in problems (\ref{abstract DA}) and (\ref{surrogate DA}) can result in problems which are not tractable using exact optimisation methods. Difficult routing problems such as the orienteering problem as utilised in the experiments as a reconnaissance decision-making model are typically solved using heuristics or constraint programming. We therefore need a solver agnostic DO method which can act as an intermediate layer. Solver agnostic methods typically rely on perturbed optimisers and stochastic interpolation to generate gradient approximations. The integrality of $\textbf{s}$ can not be relaxed as it has symbolic meaning in subsequent computation. This prevents the use of any perturbed optimisers which average over perturbed solutions. The differentiable optimisation also needs to be able to compute the gradient $\frac{\partial \textbf{s}^*}{\partial \bm{\pi}}$ directly so it can function as a layer. A popular differentiable optimisation method which fulfils all three criteria is the method proposed by \citet{poganvcic2020differentiation}. It computes the solution for a single perturbation of the cost function informed by the gradient $\frac{\partial \mathrm{loss_{DFL}}}{\partial \textbf{s}^*}$ and estimates the gradient as the difference between the perturbed solution and the initial solution. This can be interpreted as the direction of the linear interpolation between the two solutions.

For the DO gradient $\frac{\partial \textbf{x}^*}{\partial \hat{\theta}}$ the constraints that it should function as a layer and that integrality needs to be preserved in consequent computation no longer apply. Based on empirical results \citep{tang2024pyepo} we select the perturbed Fenchel-Young loss (PFYL) proposed by \citet{berthet2020learning}. PFYL is a stochastic perturbation approach which uses a theoretically sound loss function based on the Fenchel duality. 

\subsection{Optimisation Strategies}
A natural ablation is that the model parameters $W$ can be pretrained. Jointly training the DA surrogate parameters $\bm{\pi}$ and the prediction model from scratch risks early confinement to local optima basins. Not biasing the masking with a narrow range of choices for $s^*$ initially enables the model to learn from a much more diverse context. The proportion of random masking in the pretraining is set to roughly match the proportion of the masking as a result of DA constraints. We design a set of optimisation strategy ablations with the goal of improving generalisation.
\paragraph{Random search over pretrained}
One benefit of a pretrained prediction model is that it enables us to reasonably score DA strategies for randomly initialised $\pi$. Since a masked-input prediction model will work reasonably well out of the box we require only evaluation on a holdout set to score a range of random initialisation strategies for $\pi$. Initialisations which perform well can be used as warm starts for other training runs as they indicate a good DA solution.
\paragraph{Learning or fixing data acquisition}
The DO approach to optimising $\pi$ is costly as it requires repeatedly solving the DA surrogate. A key determination is whether learning $\bm{\pi}$ provides a performance advantage over fixing $\bm{\pi}$ after selecting it using a heuristic procedure such as the one proposed in the previous paragraph.
\paragraph{Joint learning or fine-tuning}
We test our hypothesis that learning the whole system for scratch leads to generalisation issues due to a less diverse training set for the prediction model. The ablation is whether to start with a pretrained prediction model or to randomly initialise one and optimise end-to-end from scratch.
\paragraph{PFL or DFL}
The final ablation is whether the decision quality is better if we optimise data acquisition for prediction loss or decision loss. For problems in which the data acquisition is fixed the choice is task dependent as it is determined by the sensitivity structure of the problem and the nature of the relationship between context and unknown parameters. \citet{cameron2022perils} show that DFL offers favourable performance under the presence of uncertainty and correlation between unknown parameters. To the best of our knowledge, the effect of missing values/masking on DFL performance has not been studied. In contrast, masked-input prediction models demonstrate robust performance on prediction tasks with missing values.
\section{Empirical Evaluation}
\subsection{Setting: Drone Reconnaissance For Shortest Path}
For our downstream contextual decision task, we borrow the Warcraft shortest path problem from \citet{poganvcic2020differentiation}. The CDT is to optimally traverse a grid based on an aerial image of it. The first step is to learn a computer vision model which predicts the cost to traverse each tile on the grid and then solve a shortest path problem from the top left to the bottom right on a $k\times k$ grid based on these travel costs. 

\begin{figure*}
\centering
\begin{subfigure}[t]{0.3\textwidth}
    \centering
    \includegraphics[width=\linewidth]{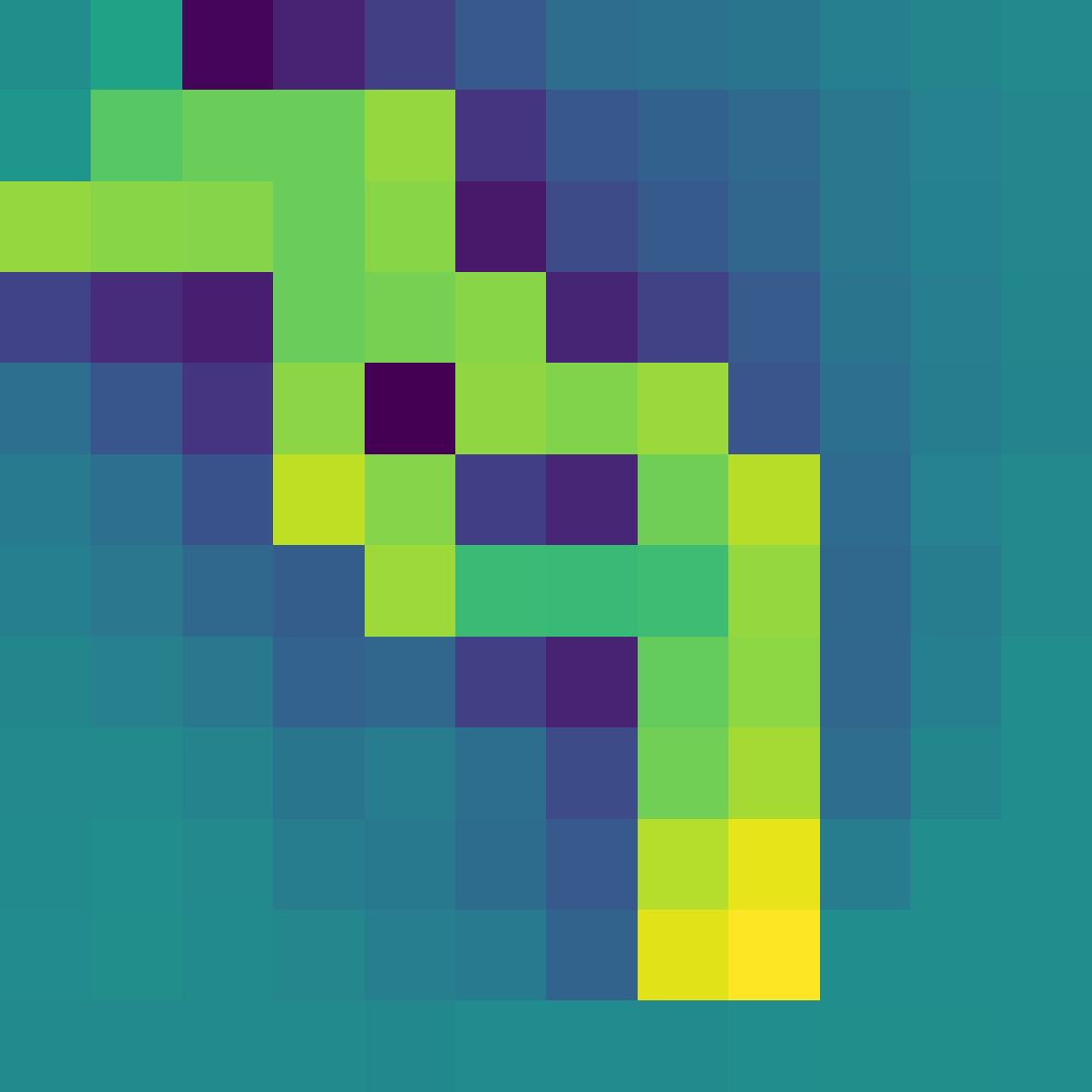}
    \caption{Learnt values of $\bm{\pi}$.}
    \label{fig:1}
\end{subfigure}
\hfill
\begin{subfigure}[t]{0.3\textwidth}
    \centering
    \includegraphics[width=\linewidth]{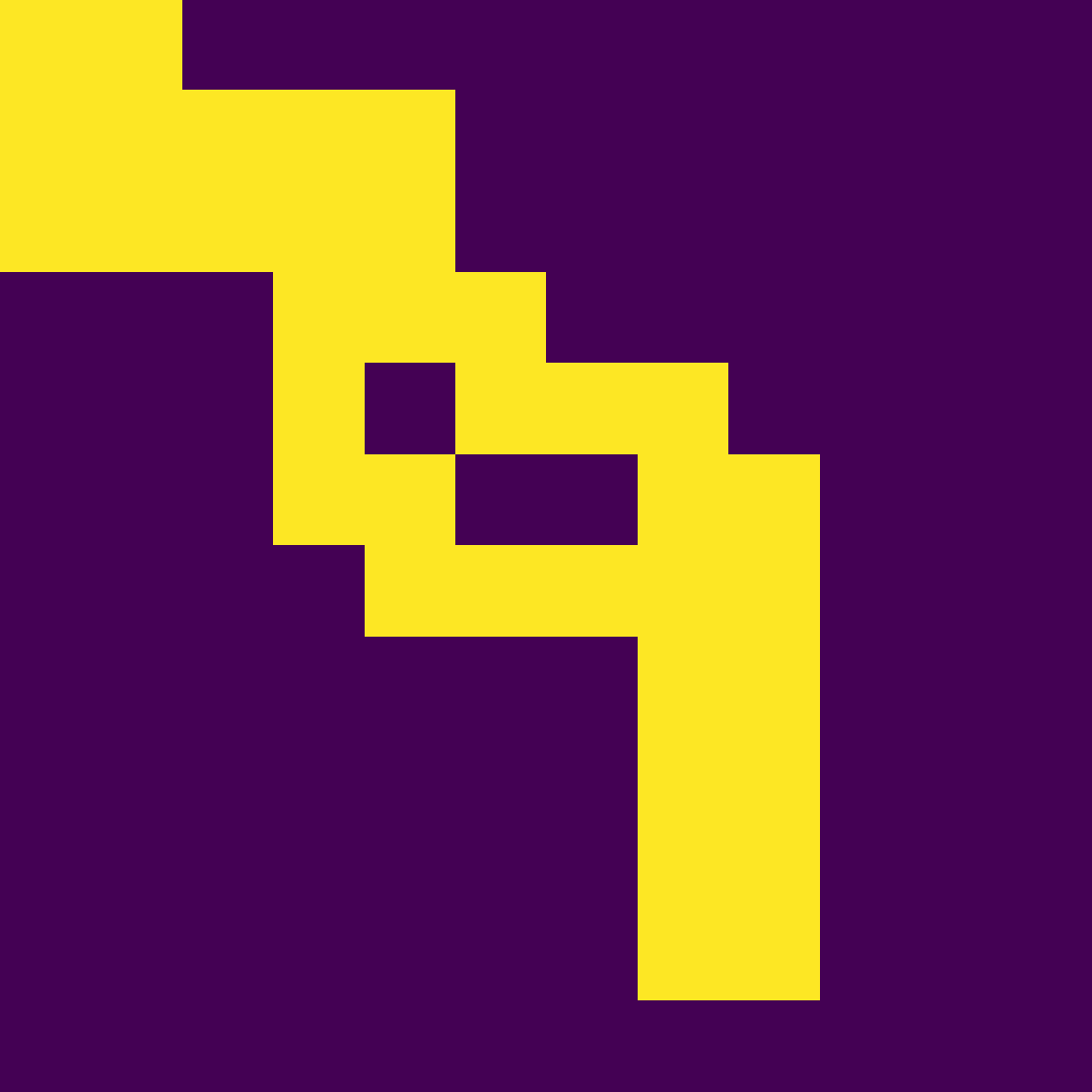}
    \caption{Solution to surrogate DA  $\textbf{s}^*$ given $\bm{\pi}$.}
    \label{fig:2}
\end{subfigure}
\hfill
\begin{subfigure}[t]{0.3\textwidth}
    \centering
    \includegraphics[width=\linewidth]{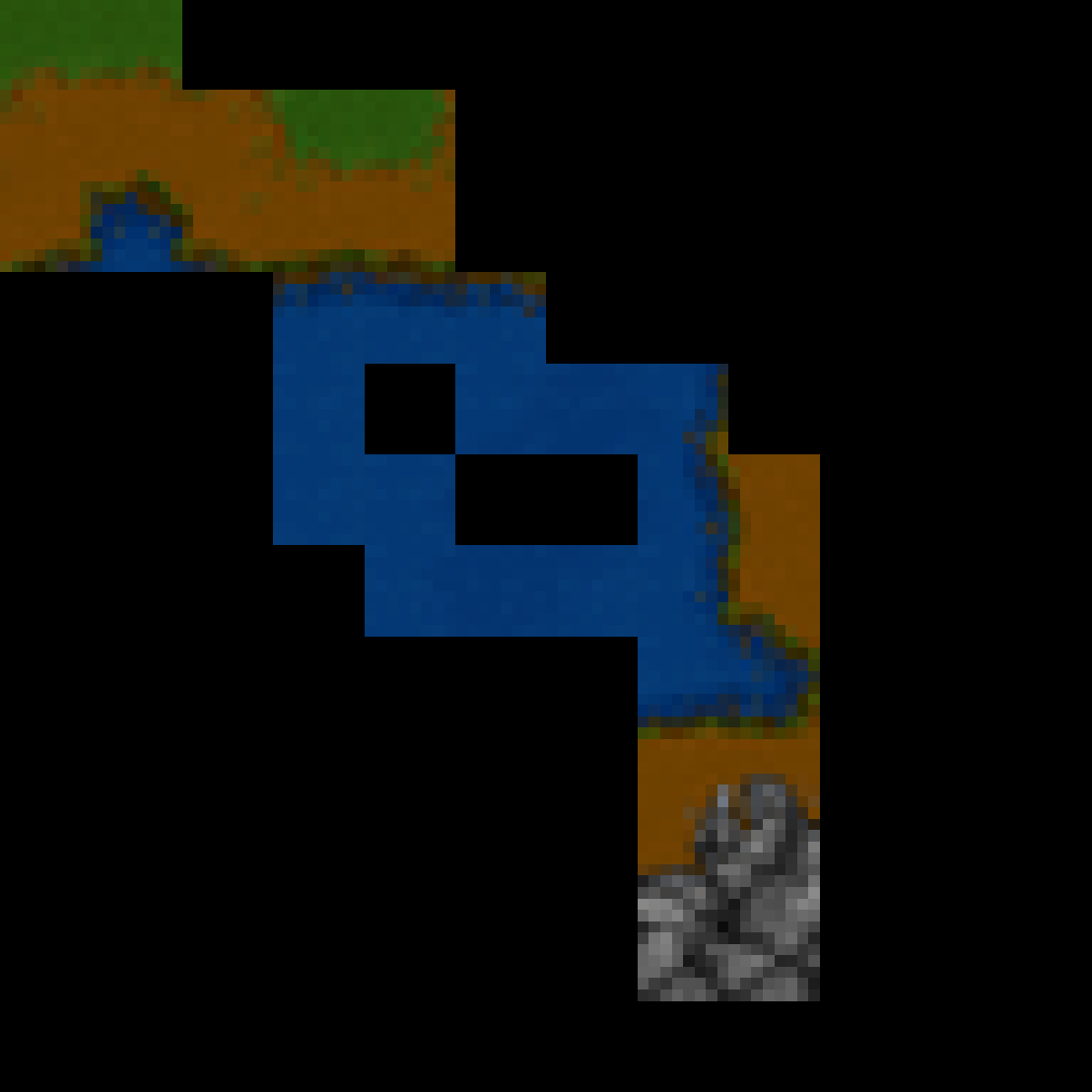}
    \caption{Observed context $\textbf{z}$ given the DA decision $\textbf{s}^*$.}
    \label{fig:3}
\end{subfigure}

\bigskip
\begin{subfigure}[t]{0.3\textwidth}
    \centering
    \includegraphics[width=\linewidth]{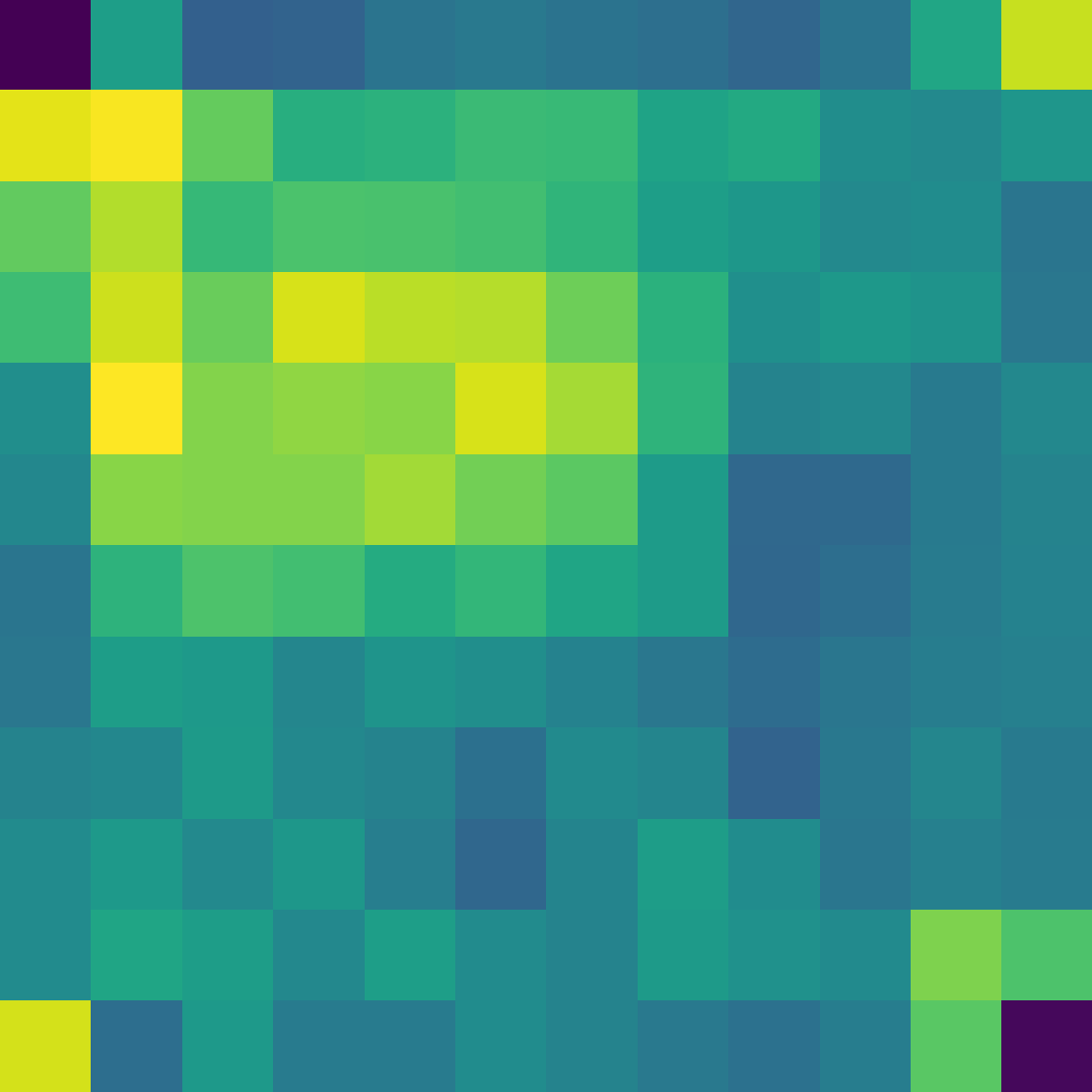}
    \caption{Predicted travel times $\hat{\theta}$ using a masked-input ViT.}
    \label{fig:4}
\end{subfigure}
\hfill
\begin{subfigure}[t]{0.3\textwidth}
    \centering
    \includegraphics[width=\linewidth]{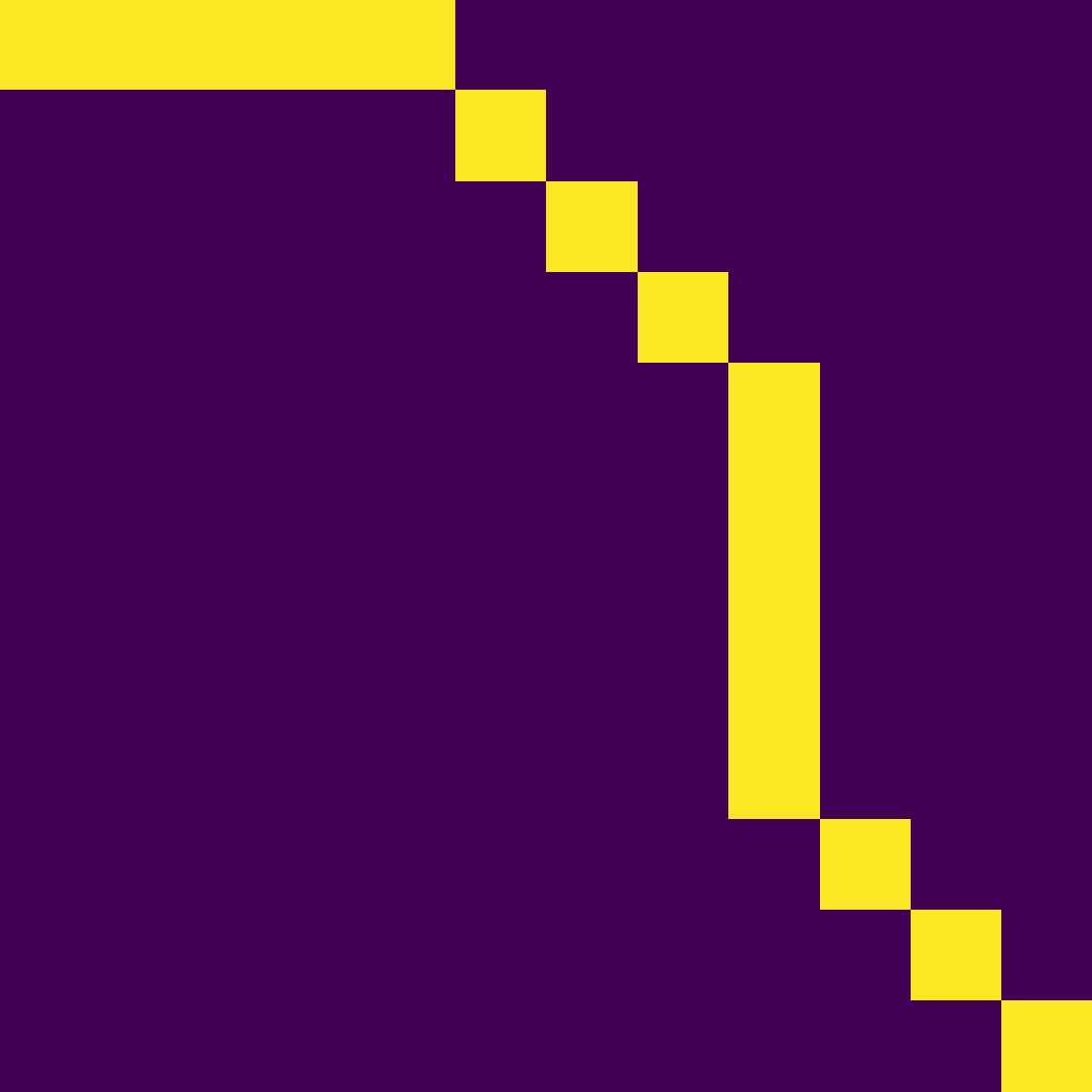}
    \caption{OPO shortest path solution $\textbf{x}^*$ obtained from $\hat{\theta}$.}
    \label{fig:5}
\end{subfigure}
\hfill
\begin{subfigure}[t]{0.3\textwidth}
    \centering
    \includegraphics[width=\linewidth]{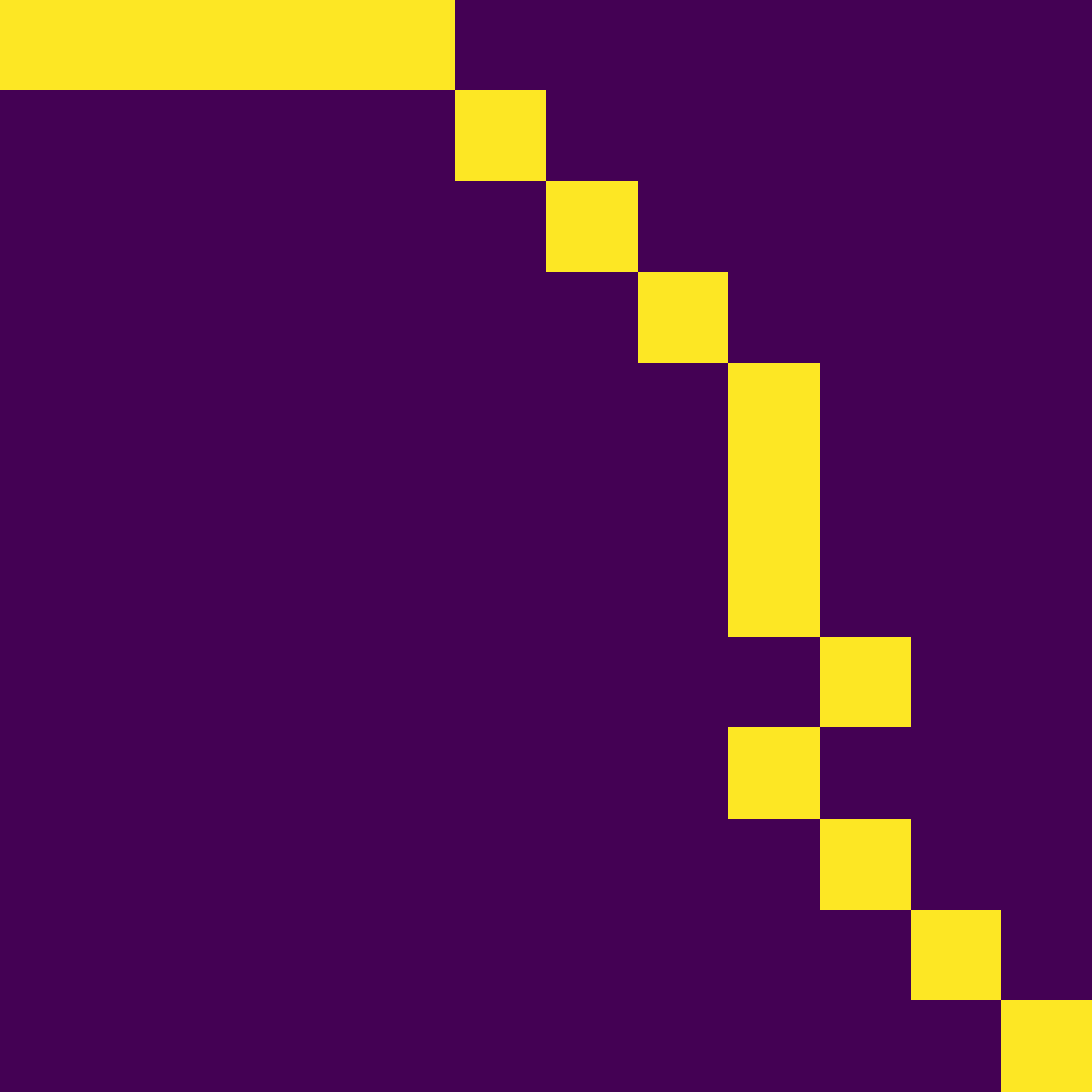}
    \caption{True shortest path solution $\textbf{x}^*(\theta)$.}
    \label{fig:6}
\end{subfigure}
\caption{OPO stages for a single instance on the Drone Reconnaissance For Shortest Path experiment from a configuration with warm-starting, learnt $\bm{\pi}$, fine-tuning, and DFL loss. Blue to yellow in \ref{fig:1} and \ref{fig:4} corresponds to low to high. Results for \ref{fig:2}, \ref{fig:5}, and \ref{fig:6} are binary: yellow corresponds to $1$ and purple to $0$.} 
\label{fig:OPO on DR}
\end{figure*}

We modify this problem by imposing an observation restriction on the test set. On the test set one has access to a single drone. The drone can only move up/down, right/left from tile to tile. It takes an $8\times 8$ RGB photo of each tile it visits. The drone has to start and finish its journey in the top left corner so it can be safely recovered. We calculate the distance between tiles as the Manhattan distance between grid coordinates reflecting the restrictions on movement directions. The drone has a range of $h < k^2$ and so not all tiles can be visited. This is equivalent to the orienteering problem, a type of routing problem in which we maximise the rewards we collect from visiting a set of nodes while constrained by travel costs. The orienteering problem is notoriously difficult to solve using mathematical programming so we utilise routing heuristics from OR-Tools \citep{ortools_routing}. The range $h$ in this case is an upper bound on how many tiles we can obtain the images for. We set $h = 36$: at most we can visit a quarter of the nodes. OPO aims to learn a set of rewards assigned to nodes such that this reconnaissance path optimises performance on the downstream shortest path problem. The solution to the surrogate DA problem is fed into our modified masked-input ViT which predicts node travel times. The predicted node travel times $\hat{\theta}$ are used to parametrise a linear programming formulation for the shortest path problem. To illustrate OPO's computation process, we present the outputs at different stages of OPO in figure \ref{fig:OPO on DR}. The code used for experiments is available  \href{https://github.com/EgoPer/OPO}{here}.

\subsection{Configuration}
We pretrain two masked-input ViTs one with a decision loss objective and the other with a prediction loss objective. We configure the ViTs such that they randomly mask three quarters of their inputs during training. We set $d_m=32$. We use the same holdout splits as \citet{poganvcic2020differentiation} which contain $10000$ training examples, and $1000$ validation and test examples each. We train both until an early stopping criterion of no validation set performance improvement in 10 epochs is met. For prediction model parameters we use a learning rate of $3\text{e-}4$ with an Adam optimiser throughout all experiments.

We initialise $\bm{\pi}$ as independent draws from a uniform distribution between $0.495$ and $0.505$. The ortools solver for the orienteering problem requires integral non-negative values. In the optimiser we multiply the values of $\bm{\pi}$ by $1\text{e}3$, round them, and pass them through a ReLU layer. We also add the value of $h$ to the reward at each node so the travel costs do not influence the reward collection process. Note that in this case the optimiser computations do not affect the backpropagation. The initialisation was selected based on informal observations on the scale of gradients and to allow for a diverse set of reward values.
\subsection{Results}
\begin{figure*}
\centering
\begin{subfigure}[t]{0.49\textwidth}
    \centering
    \includegraphics[width=\linewidth]{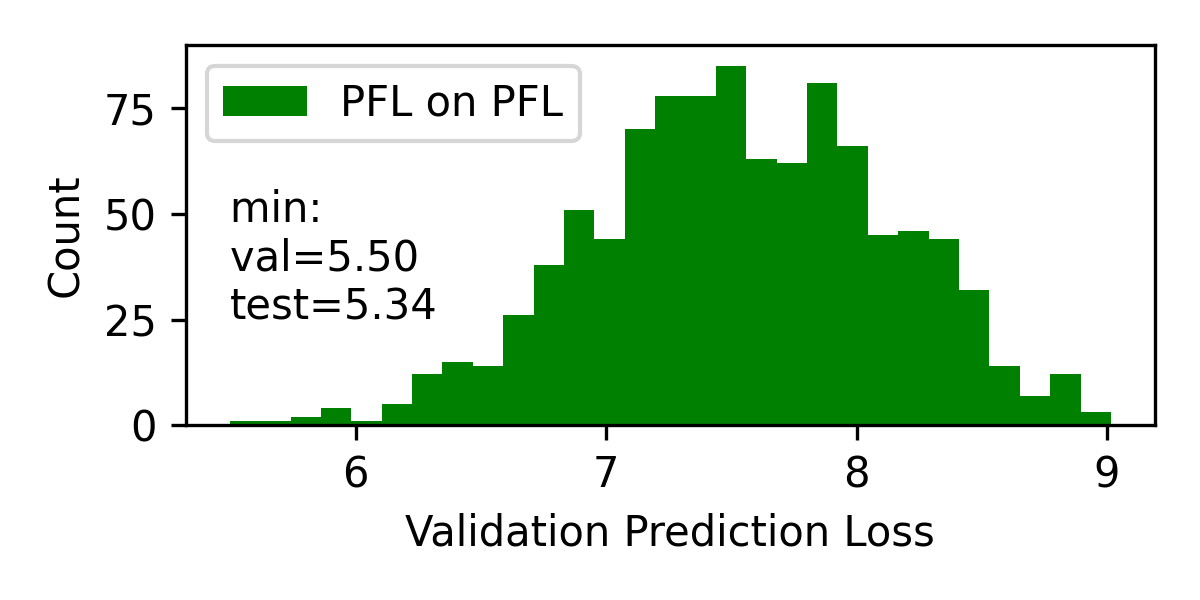}
    \caption{PFL task.}
    \label{fig: rs pfl}
\end{subfigure}
\hfill
\begin{subfigure}[t]{0.49\textwidth}
    \centering
    \includegraphics[width=\linewidth]{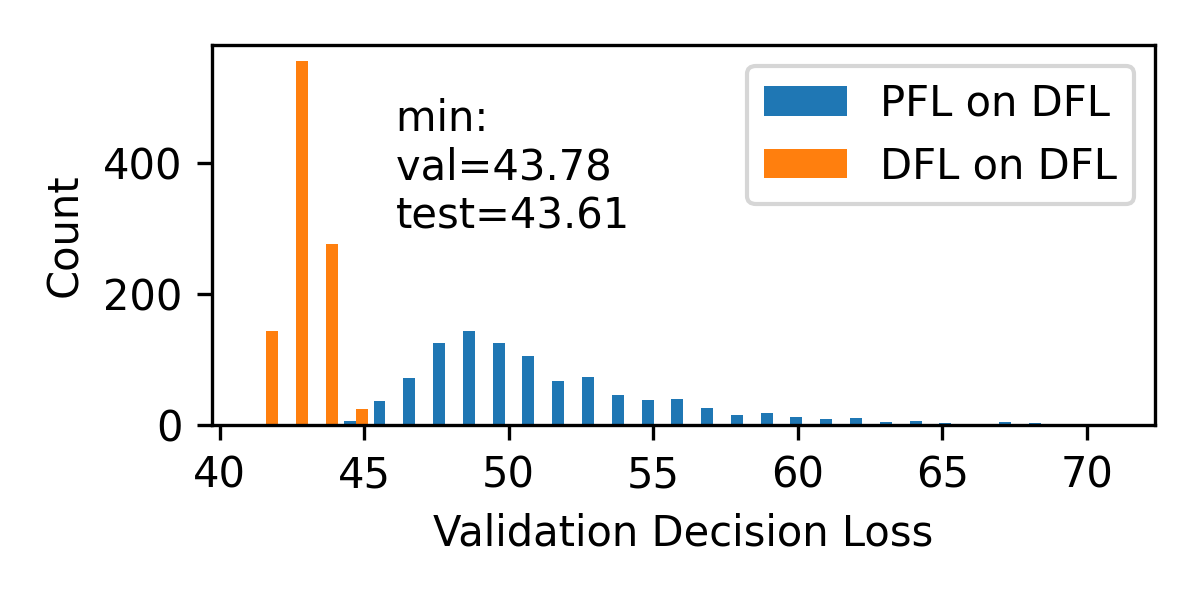}
    \caption{DFL task.}
    \label{fig: rs dfl}
\end{subfigure}
\caption{Distribution of mean validation set results for the random search strategy. We note the best performing initialisation. We do not list results of DFL trained on PFL as they are obviously not competitive.} 
\label{fig: RS}
\end{figure*}

\paragraph{Random Search} We initialise $\bm{\pi}$ with $1000$ seeds and compute the validation scores for each model and loss criterion. We present the distribution of performance for training/task combinations in figure \ref{fig: RS}. For comparison note that the mean objective value given a-priori knowledge of the CDT parameters is $29.64$ on the validation set and $29.71$ on the test set. None of the initialisations result in a DA strategy which makes a PFL pretrained model competitive on the decision task. The spread of the results suggests that some DA strategies are substantially better than others thus warranting optimisation on this problem. Notably, the variance of mean objectives for the PFL trained model when evaluated for decision quality is substantially higher than any other combination. This indicates that downstream decisions are highly sensitive to missing data if the prediction model is not trained end-to-end or fine-tuned to a specific data acquisition strategy. 

\paragraph{Ablations} We select the top four performing initialisations across both criteria in addition to a fifth initialisation with a roughly average performance to warm-start all other optimisation approaches. We perform one training run across all ablation combinations. In the initial training run we do not perform any hyperparameter optimisation and set the learning rate for $\bm{\pi}$ at $1\text{e-}3$. We present the results of this experiment in Table \ref{table: initial results}.

\begin{table*}[t]
\caption{Initial training run results with training across five initialisations. Each combination reflects runs for five initialisations. The "best" reflects selection amongst initialisations based on the validation score on that specific task. Under the $W$ columns FT stands for the fine-tuning ablation and JO stands for the joint learning ablation.}
\label{table: initial results}
\vskip 0.15in
\begin{center}
\begin{small}
\begin{sc}
\begin{tabular}{l|l|c|c|c|c|c|c|c|c|c|c|c|c}
\toprule
& & \multicolumn{8}{c|}{pfl training} & \multicolumn{4}{c}{dfl training} \\
\cline{3-14}
 &  & 
 \multicolumn{4}{c|}{$\mathrm{loss_{PFL}}$} &
 \multicolumn{4}{c|}{$\mathrm{loss_{DFL}}$} &
 \multicolumn{4}{c}{$\mathrm{loss_{DFL}}$} \\
\cline{3-14}
$\bm{\pi}$ & $W$ & mean & std & min & best & mean & std & min & best & mean & std & min & best \\
\midrule
\multirow[c]{2}{*}{fixed} & FT 
& 4.37 & 0.56 & 3.92 & \textbf{3.92} & 38.81 & 0.97 & 37.40 & 39.66
& 39.44 & 0.75 & 38.69 & 38.74 \\
& JO 
& 4.59 & 0.54 & 4.14 & 4.14 & 39.74 & 0.29 & 39.32 & 40.11 
& 39.27 & 0.72 & 38.29 & \underline{38.29}
 \\
\midrule
\multirow[c]{2}{*}{learnt} & FT 
& 5.01 & 0.74 & 4.32 & 4.32 & 39.99 & 1.98 & 37.42 & \textbf{37.42}
& 40.29 & 1.15 & 38.98 & 38.98 \\
& JO 
& 5.37 & 0.52 & 4.85 & 4.93 & 40.08 & 0.58 & 39.45 & 39.45  
& 43.30 & 1.89 & 39.93 & 39.93 \\
\bottomrule
\end{tabular}
\end{sc}
\end{small}
\end{center}
\vskip -0.1in
\end{table*}

We see a substantial improvement in the performance of masked-input prediction models as a result of fine-tuning which is consistent with expectations and in some of the cases of joint optimisation over random search for DA strategies. The best overall result is achieved by a learnt DA strategy presenting evidence in favour OPO over random search and fine-tuning. PFL outperforms DFL on this task for a learnt $\bm{\pi}$ which is consistent with empirical results on this specific task \citep{tang2024pyepo}. This is not true for fixed $\bm{\pi}$ which may reflect the bias specific to the observed context. We hypothesise that the inconsistency between prediction and decision loss in PFL training reveals a degree of overfitting in favour of the prediction task given a specific DA strategy. In contrast learning the DA means that the masked-input model will see a more diverse set of contexts initially before settling on a DA strategy an fine-tuning. This is further evidenced by validation selection being more consistent with the best overall model in the learnt setting.

When training with a DFL objective the best performing combination is joint learning with a fixed DA strategy indicating the DFL-pretrained masked-input model may have overfitting issues for the warm-start DA strategies. Fine-tuning appears to have a larger effect on learnt DA in the PFL training setting both in terms of the best result and performance variance. This and the fact that some training runs saw no performance improvements suggests that there is a degree of instability in the optimisation of $\bm{\pi}$ especially when fine-tuning. To see if we can stabilise the training in the fine-tuning ablation we tune the learning rate for $\bm{\pi}$ on the initialisation with the best validation score across experiments in both the PFL and DFL setting. We try eight learning rates for $\bm{\pi}$ with roughly logarithmic spacing. The results are presented in Table \ref{table: tuned pi}.

\begin{table*}[t]
\caption{Results for eight learning rates 
for $\bm{\pi}$ on the best performing initialisation from the initial experiment.}
\label{table: tuned pi}
\vskip 0.15in
\begin{center}
\begin{small}
\begin{sc}
\begin{tabular}{l|l|c|c|c|c|c|c|c|c|c|c|c|c}
\toprule
& & \multicolumn{8}{c|}{pfl training} & \multicolumn{4}{c}{dfl training} \\
\cline{3-14}
 &  & 
 \multicolumn{4}{c|}{$\mathrm{loss_{PFL}}$} &
 \multicolumn{4}{c|}{$\mathrm{loss_{DFL}}$} &
 \multicolumn{4}{c}{$\mathrm{loss_{DFL}}$} \\
\cline{3-14}
$\bm{\pi}$ & $W$ & mean & std & min & best & mean & std & min & best & mean & std & min & best \\
\midrule
\multirow[c]{1}{*}{learnt} & FT 
& 4.75 & 0.42 & 4.34 & 4.34 & 38.55 & 0.94 & 36.84 & \textbf{36.84}
& 38.78 & 0.55 & 38.08 & \underline{38.29} \\
\bottomrule
\end{tabular}
\end{sc}
\end{small}
\end{center}
\vskip -0.1in
\end{table*}

Tuning the learning rate results in an improvement in performance under both training modalities. The best overall performance corresponds to a learning rate of $7\text{e-}4$ for the DA surrogate. In terms of relative regret compared to having perfect information the decision performance is $0.469$ under random search, to $0.289$ under the best performing fixed $\bm{\pi}$ instance, to $0.260$ for the best performing learnt $\bm{\pi}$ before tuning, and to $0.240$ after tuning. The relative regret of the best OPO model is a substantial $17\%$ lower than the best non-learnt method.

\paragraph{Scaling} The additional training time required for OPO depends on the choice of DO method and optimisation approach used to solve the surrogate DA task. Our implementation uses a heuristic solver with a ten second time budget. Each batch requires two solutions to the surrogate DA task, once in the forward and once in the backward pass. The additional twenty seconds per batch are not prohibitive but such short heuristic run times would not yield good solutions on larger surrogate DA tasks. As with DFL, solver calls are the computationally expensive element. The development of solver-free approaches which work as a layer would be of immense benefit.

\section{Related Work}
\paragraph{Existing systematic DA}
Optimal sensor placement concerns itself with how to optimally place a limited number of sensors across a geography to optimise some coverage objective. A popular family of solutions is the use of Gaussian Processes to build a surrogate model which is then queried for objectives such as mutual information \citep{krause2008near} commonly used due to its submodularity. Sensor placement has trivial cardinality constraints. In contrast informative path planning (IPP) involves routing constraints which demand more complex optimisation approaches such as branch-and-bound \citep{binney2012branch}. IPP is a mature field exemplified by efficient algorithms for online learning in unknown environments \citep{schmid2020efficient} using sampling techniques such as rapidly exploring random trees. The online approaches for IPP may contain insight on how to extend OPO to environments without the full observation assumption. The key difference between OPO and existing approaches such as IPP is that OPO explicitly optimises for downstream objectives instead of relying on statistical proxies thus avoiding the risk of objective misalignment.

\paragraph{Surrogate Modelling}
Surrogate modelling utilising differentiable optimisation (DO) has been employed to solve hard optimisation problems. \citet{ferber2023surco} propose SurCO, a surrogate mixed-integer programming (MIP) model with a learnable linear objective to approximate the optimal solutions of non-linear MIP problems with equivalent constraints. For a fixed prediction model, and MIP constraints in the DA problem OPO can be interpreted as applying SurCo to an implicit non-linear objective.
\section{Conclusion}
OPO is the first model to tackle the data acquisition problem end-to-end enabling a greater alignment between DA decisions and the true downstream objective of a decision-maker. OPO is easily extendable to contextual DA settings in which some context is available before the DA decision has to be made by prepending a prediction layer. We believe dynamic learning rates could be leveraged to further improve performance due to the observed sensitivity of the surrogate DA problem. To the best of our knowledge this is the first instance of non-sequential optimisation layers in a neural architecture. Our experiments on the drone reconnaissance problem provide evidence in support of broader end-to-end modelling in multilevel optimisation problems. 






\bibliography{references}
\bibliographystyle{icml2025}

\end{document}